\documentclass[11pt]{amsart}
\usepackage{amssymb}

\input{prepictex.tex}
\input{pictex.tex}
\input{postpictex.tex}

\begin{document} 

\title[Graph algebras]{Stable Rank of Graph Algebras. \\ 
       Type $\mathbf{I}$ Graph Algebras and their Limits} 

\author[Deicke]{Klaus Deicke} 
\address{Fachbereich 17/Mathematik \\ 
Universit\"at Paderborn \\ 
D--33095 Paderborn \\ 
Germany} 
\email{deicke@uni-paderborn.de} 

\author[Hong]{Jeong Hee Hong} 
\address{Department of Applied Mathematics \\ 
Korea Maritime University \\ 
Busan 606--791 \\ 
South Korea} 
\email{hongjh@hanara.kmaritime.ac.kr} 

\author[Szyma\'{n}ski]{Wojciech Szyma\'{n}ski} 
\address{Mathematics \\
The University of Newcastle \\
NSW 2308 \\
Australia}
\email{wojciech@frey.newcastle.edu.au} 


\thanks{This research was partially supported by a grant from the 
exchange program of the Australian Academy of Science and Korea Science 
and Engineering Foundation.} 

\date{31 July, 2002} 

\def\T{\mathbb T}
\def\Z{\mathbb Z}
\def\R{\mathbb R}
\def\C{\mathbb C} 
\def\D{\mathbb D} 
\def\N{\mathbb N}

\def\B{\mathcal B} 
\def\H{\mathcal H}
\def\K{\mathcal K}
\def\F{\mathcal F}
\def\U{\mathcal U}
\def\E{\mathcal E}
\def\O{\mathcal O}
\def\J{\mathcal J} 

\def\Xbad{X^{\text{{\rm fin}}}_{\infty}} 
\def\OMbad{\Omega(M)^{\text{{\rm fin}}}_{\infty}}
\def\Ovbad{\Omega(v)^{\text{{\rm fin}}}_{\infty}}
\def\Oabad{\Omega(\alpha)^{\text{{\rm fin}}}_{\infty}}

\def\Aut{\operatorname{Aut}}
\def\Prim{\operatorname{Prim}}
\def\id{\text{id}}

\newtheorem{defi}{Definition}[section]
\newtheorem{coro}[defi]{Corollary}
\newtheorem{lemm}[defi]{Lemma}
\newtheorem{prop}[defi]{Proposition}
\newtheorem{theo}[defi]{Theorem}
\newtheorem{exam}[defi]{Example}
\newtheorem{rema}[defi]{Remark}

\begin{abstract}
For an arbitrary countable directed graph $E$ we show that the 
only possible values of the stable rank of the 
associated Cuntz-Krieger algebra $C^*(E)$ are $1$, 
$2$ or $\infty$. Explicit criteria for each of these 
three cases are given. We characterize graph algebras of type $I$, 
and graph algebras which are inductive limits of $C^*$-algebras 
of type $I$. We also show that a gauge-invariant ideal of a 
graph algebra is itself isomorphic to a graph algebra. 
\end{abstract} 

\maketitle

\setcounter{section}{-1}

\section{Introduction} 

In the last few years a great progress has been made in the 
research of graph $C^*$-algebras. In spite of the fact that this 
class contains objects so diverse as classical Cuntz-Krieger algebras \cite{ck}, 
$AF$-algebras, and a variety of type $I$ algebras, uniform methods have 
been developed for dealing with the algebras associated to arbitrary 
countable directed graphs. The purpose of this article is to make further 
progress toward understanding of their structure. In our main result, Theorem 
\ref{main}, we determine the (topological) stable rank of the algebra 
$C^*(E)$ associated to an arbitrary countable graph $E$. Indeed, we 
show that it can take only values $1$, $2$ and $\infty$, and give 
criteria for each of these three cases. Previously only the case of 
stable rank one was satisfactorily settled in \cite[Proposition 5.5]{rs} 
(see also \cite[Theorem 3.3]{jps} for the special case of row-finite graphs). 

As a key step in our argument we show that if a graph algebra does not admit 
unital, purely infinite and simple quotients then it has a stable ideal 
whose quotient is a limit of type $I$ algebras (cf. Lemma \ref{stableideal}). 
To carry out this idea we give characterizations of graph algebras of type 
$I$ (cf. Theorem \ref{typeI}), and graph algebras which are isomorphic to 
inductive limits of $C^*$-algebras of type $I$ (cf. Theorem \ref{limitoftypeI}). 
The latter ones constitute a very interesting subclass of graph algebras. They 
are easily identifiable as the algebras of graphs `with isolated loops', in the 
sense to be made precise in \S1. They appear to be exactly on the opposite 
side of the whole spectrum of graph algebras to the purely infinite ones 
(in the sense of \cite{kr}) discussed in \cite{hs4}. It seems logical that 
future classification of graph algebras should begin with these two extreme 
subclasses. 

The results of this paper depend heavily on the vast machinery of graph 
algebras developed up-to-date. The most important for us was the 
classification of their ideals, recently completed in \cite{bhrs} and 
\cite{hs3}, and the method of approximating algebras of infinite graphs 
with the algebras of finite graphs, introduced in \cite{rs}. Among others, 
we have also used the general criterion of injectivity of homomorphisms of 
graph algebras, proved in \cite{s2}, and the characterization of 
stable graph algebras, given in \cite{t}. Of course, our computation 
of stable rank depends crucially on the numerous results of the Rieffel's 
original paper \cite{r}. 

One technical lemma of independent interest proved in the present article 
says that every gauge-invariant ideal of a graph algebra is itself isomorphic 
to a graph algebra (cf. Lemma \ref{Jideal}). Combining this with 
Proposition 3.4 and Corollary 3.5 of \cite{bhrs} we see that the class 
of graph algebras is closed under taking gauge-invariant ideals and 
their quotients. This is yet another evidence of the well-known but 
perhaps not so well understood fact that generalized Cuntz-Krieger 
algebras and their gauge actions go together. 

\section{Preliminaries on graph algebras} 

We recall the definition of the $C^*$-algebra 
corresponding to a directed graph \cite{flr}. 
Let $E=(E^0,E^1,r,s)$ be a directed graph with 
countably many vertices $E^0$ and  
edges $E^1$, and range and source functions 
$r,s:E^1\rightarrow E^0$, respectively. The algebra 
$C^*(E)$ associated to the graph $E$ is defined as the universal 
$C^*$-algebra generated by families of projections 
$\{P_v:v\in E^0\}$ and partial isometries $\{S_e:e\in E^1\}$, 
subject to the following relations.  
\begin{description} 
\item[(GA1)] $P_vP_w=0$ for $v,w\in E^0$, $v\neq w$. 
\item[(GA2)] $S_e^*S_f=0$ for $e,f\in E^1$, $e\neq f$. 
\item[(GA3)] $S_e^*S_e=P_{r(e)}$ for $e\in E^1$.  
\item[(GA4)] $S_e S_e^*\leq P_{s(e)}$ for $e\in E^1$.  
\item[(GA5)] $P_v=\displaystyle{\sum_{e\in E^1,
s(e)=v}}S_e S_e^*$ for $v\in E^0$ such that $0<|s^{-1}(v)|<\infty$. 
\end{description} 
Universality in this definition means that if $\{Q_v:v\in E^0\}$ and 
$\{T_e:e\in E^1\}$ are families of projections and partial 
isometries, respectively, satisfying conditions (GA1--GA5), then 
there exists a $C^*$-algebra homomorphism from $C^*(E)$ to the 
$C^*$-algebra generated by $\{Q_v:v\in E^0\}$ and 
$\{T_e:e\in E^1\}$ such that $P_v\mapsto Q_v$ and 
$S_e\mapsto T_e$ for $v\in E^0$, $e\in E^1$. Universality implies 
existence of the gauge action $\gamma:\T\rightarrow\Aut(C^*(E))$  
such that $\gamma_z(P_v)=P_v$ and $\gamma_z(S_e)=zS_e$ for all 
$v\in E^0$, $e\in E^1$, $z\in\T$. 

More than once in this article we will use the fact that the $C^*$-algebra 
of an arbitrary infinite graph may be written as the increasing limit  
of its $C^*$-subalgebras which are themselves isomorphic to the algebras 
of finite graphs. This fact was proved in \cite[Lemma 1.2 and \S5.1]{rs}. 
However, because of the importance of this approximation to our approach 
we want to recall the details of that construction. Our present approach 
differs slightly from that of \cite[Definition 1.1]{rs} by the inclusion 
of both edges and vertices in the generating subfamily, as well as by 
the possible presence of sinks in the original graph. 

So let $E$ be an arbitrary countable directed graph. Let $G\subseteq 
E^0\cup E^1$ be a finite set. We denote $G^0=G\cap E^0$ and $G^1=G\cap E^1$,  
assuming that $G^0$ contains $r(G^1)$. We define a finite graph $E_G$ 
with vertices $E_G^0=(E_G)^0$ and edges $E_G^1=(E_G)^1$, as follows: 
\begin{align*} 
E_G^0 &:= G^1\cup\{v\in G^0: \text{ either } s^{-1}(v)=\emptyset 
          \text{ or } v\in s(E^1\setminus G^1)\}, \\ 
E_G^1 &:= \{(e,f)\in G^1\times E_G^0:r(e)=s(f)\}, 
\end{align*} 
with $s(e,f)=e$ and $r(e,f)=f$. Note that each vertex $v\in E_G^0
\setminus G^1$ is a sink in $E_G$. In the following Lemma \ref{approximation} 
we denote the canonical Cuntz-Krieger family generating $C^*(E_G)$ by 
$\{T_*,Q_*\}$. The proof of the lemma is very similar to that of 
\cite[Lemma 1.2]{rs} and hence it is omitted. 

\begin{lemm}\label{approximation} 
Let $E$ be a directed graph and let $G\subseteq E^0\cup E^1$ be finite. 
Then there exists an isomorphism $\phi_G$ from $C^*(E_G)$ onto the 
$C^*$-subalgebra of $C^*(E)$ generated by $\{S_e:e\in G^1\}$ and $\{P_v
:v\in G^0\}$ such that 
\begin{align*} 
\phi_G(Q_e) &= S_e S_e^*, \\ 
\phi_G(Q_v) &= P_v-\sum_{f\in G^1,s(f)=v}S_f S_f^*, \\ 
\phi_G(T_{(g,h)}) &= S_g S_h S_h^*, \\ 
\phi_G(T_{(k,w)}) &= S_k\left(P_w-\sum_{f\in G^1,s(f)=w}S_f S_f^*\right), \\ 
\end{align*} 
for $e,g,h,k\in G^1$, $v,w\in E_G^0\setminus G^1$, such that $r(g)=s(h)$, 
$r(k)=w$. 
\end{lemm} 

Choosing an exhaustive increasing sequence of finite subsets $G_n$ of $E^0\cup 
E^1$ and applying Lemma \ref{approximation}, we see that $C^*(E)$ is an 
increasing limit of its $C^*$-subalgebras isomorphic to the algebras 
$C^*(E_{G_n})$, with each $E_{G_n}$ a finite graph. 

\vspace{2mm} 
It is useful to notice that the graph $E_G$ preserves certain important 
properties of the original graph $E$ (cf. \cite[Lemma 1.3]{rs}). In 
particular, $E_G$ preserves the property of `having isolated 
loops', which is investigated in \S2. We say that a graph $E$ has 
{\em isolated loops} if whenever $(a_1,\ldots,a_k)$ and $(b_1,\ldots,b_n)$ 
are loops in $E$ such that $s(a_i)=s(b_j)$ then $a_i=b_j$. At first we 
note the following. Every loop in $E_G$ has the form $\alpha=
((e_1,e_2),(e_2,e_3),\ldots,(e_k,e_1))$, with $e_j\in E^1$ such that 
$r(e_j)=s(e_{j+1})$ (and $r(e_k)=s(e_1)$), by the definition of $E_G^1$. 
Hence $(e_1,\ldots,e_k)$ is the loop in $E$ corresponding to the loop 
$\alpha$ in $E_G$. 

\begin{lemm}\label{preserve} 
Let $E$ be a directed graph and let $G\subseteq E^0\cup E^1$ be finite. If 
$E$ has isolated loops then so has $E_G$. 
\end{lemm} 
\begin{proof} 
Suppose that $((e_1,e_2),(e_2,e_3),\ldots,(e_k,e_1))$ and $((f_1,f_2),
(f_2,f_3),\ldots,(f_n,f_1))$ are two loops in $E_G$ such that 
$s((e_1,e_2))=s((f_1,f_2))$ but $(e_1,e_2)\neq(f_1,f_2)$. This 
implies $e_1=f_1$ and $e_2\neq f_2$. Consequently, $(e_1,e_2,\ldots,e_k)$ 
and $(f_1,f_2,\ldots,f_n)$ are two loops in $E$ such that $s(e_2)=s(f_2)$ 
but $e_2\neq f_2$. 
\end{proof} 

In order to determine the stable rank of the graph 
$C^*$-algebras we will use Lemma \ref{Jideal}, below, which is of 
independent interest. Combined with \cite[Theorem 3.6]{bhrs}, Lemma 
\ref{Jideal} says that every gauge-invariant ideal of a graph algebra is 
itself isomorphic to a graph algebra. In order to define the graph 
corresponding to a gauge-invariant ideal we need the following notation. 
Let $E$ be a directed graph, let $X\neq\emptyset$ be a hereditary and 
saturated subset of $E^0$ (cf. \cite{bhrs}) and let $B\subseteq\Xbad$. 
(Recall from \cite{bhrs} that $\Xbad$ is the subset of $E^0\setminus X$ 
consisting of all those vertices $v$ such that $|s^{-1}(v)|=\infty$ and 
$0<|s^{-1}(v)\cap r^{-1}(E^0\setminus X)|<\infty$.) We 
denote by $\widetilde{F}_E(X,B)$ the collection of all finite paths 
$\alpha=(a_1,\ldots,a_{|\alpha|})$ of positive length such that 
$s(\alpha)\in E^0\setminus X$, $r(\alpha)\in X\cup B$, and 
$r(a_j)\not\in X\cup B$ for $j<|\alpha|$. We define 
$$ F_E(X,B)=\widetilde{F}_E(X,B)\setminus
   \{e\in E^1:s(e)\in B \text{ and } r(e)\in X\}. $$ 
We denote by $\overline{F}_E(X,B)$ another copy of $F_E(X,B)$. If 
$\alpha\in F_E(X,B)$ then we write $\overline{\alpha}$ for the copy of 
$\alpha$ in $\overline{F}_E(X,B)$. 

\begin{defi}\label{idealgraph} 
{\rm Let $E$ be a directed graph. Let $X\neq\emptyset$ be a hereditary and 
saturated subset of $E^0$, and let $B\subseteq\Xbad$. We define a directed 
graph $_XE_B$, as follows: 
\begin{align*} 
( {_XE_B})^0= {_XE_B^0} & := X\cup B\cup F_E(X,B), \\  
( {_XE_B})^1= {_XE_B^1} & := \{e\in E^1:s(e)\in X\} \cup 
    \{e\in E^1:s(e)\in B \text{ and } r(e)\in X\} \cup\overline{F}_E(X,B), 
\end{align*}  
with $s(\overline{\alpha})=\alpha$ and $r(\overline{\alpha})=r(\alpha)$ for 
$\alpha\in F_E(X,B)$, and the source and range as in $E$ for the other 
edges of $_XE_B^1$.}  
\end{defi} 

\begin{exam}\label{example1} 
{\rm Let $E$ be the following graph (where the symbol $(\infty)$ indicates 
that there are infinitely many edges from $b$ to $x_1$): 

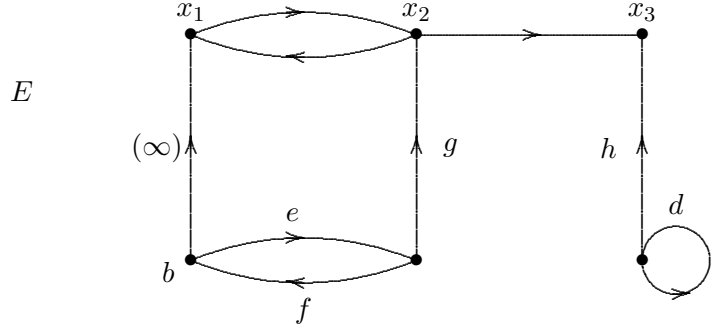
\[ \beginpicture
\setcoordinatesystem units <1.5cm,1.5cm>
\setplotarea x from -4 to 1, y from -0.7 to 0.7 

\put {$E$} at -4.5 0.5
\put {$\bullet$} at -3 -1
\put {$\bullet$} at -3 1
\put {$\bullet$} at -1 -1
\put {$\bullet$} at -1 1
\put {$\bullet$} at 1 -1
\put {$\bullet$} at 1 1

\put {$e$} at -2.1 -0.6
\put {$f$} at -2 -1.45
\put {$(\infty)$} at -3.3 0
\put {$g$} at  -0.7 0
\put {$h$} at 0.7 0
\put {$d$} at 1.3 -0.5 
\put {$x_1$} [v] at -3 1.2 
\put {$x_2$} [v] at -1 1.2 
\put {$x_3$} [v] at 1 1.2 
\put {$b$} [v] at -3.2 -1.1 

\setlinear  
\plot -3 -1 -3 1 /
\plot -1 -1 -1 1 /
\plot -1 1 1 1 / 
\plot 1 -1 1 1 /
\arrow <0.235cm> [0.2,0.6] from 0 1 to 0.1 1 
\arrow <0.235cm> [0.2,0.6] from 1 0 to 1 0.1 
\arrow <0.235cm> [0.2,0.6] from  -3 0 to -3 0.1  
\arrow <0.235cm> [0.2,0.6] from -1 0 to -1 0.1 

\setquadratic 
\plot -3 1 -2 0.8 -1 1 /
\plot -3 1 -2 1.2 -1 1 /
\plot -3 -1 -2 -1.2 -1 -1 /
\plot -3 -1 -2 -0.8 -1 -1 /

\arrow <0.235cm> [0.2,0.6] from -2 0.8 to -2.15 0.815
\arrow <0.235cm> [0.2,0.6] from -2.15 1.2 to -2 1.2
\arrow <0.235cm> [0.2,0.6] from -2.15 -0.8 to -2 -0.8
\arrow <0.236cm> [0.2,0.6] from -2 -1.2 to -2.15 -1.185

\circulararc 360 degrees from 1 -1 center at 1.3 -1  
\arrow <0.236cm> [0.2,0.6] from 1.27 -1.3 to 1.42 -1.29
 
\endpicture \] 
If $X=\{x_1,x_2,x_3\}$ and $B=\{b\}$, then the graph $_XE_B$ is: 

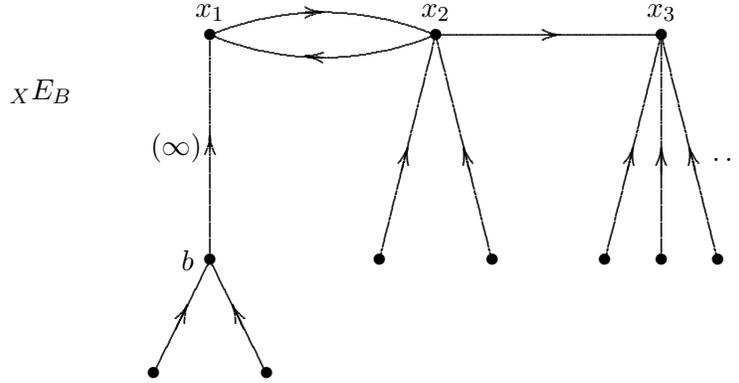
\[ \beginpicture
\setcoordinatesystem units <1.5cm,1.5cm>
\setplotarea x from -4 to 1, y from -0.7 to 0.7 

\put {$_XE_B$} at -4.5 0.5
\put {$\bullet$} at -3 -1
\put {$\bullet$} at -3 1
\put {$\bullet$} at -3.5 -2
\put {$\bullet$} at -2.5 -2 
\put {$\bullet$} at -1 1
\put {$\bullet$} at -0.5 -1
\put {$\bullet$} at -1.5 -1
\put {$\bullet$} at 1 1
\put {$\bullet$} at 1 -1
\put {$\bullet$} at 0.5 -1
\put {$\bullet$} at 1.5 -1 
\put {$(\infty)$} at -3.3 0

\setlinear  
\plot -3 -1 -3 1 /
\plot -1 1 1 1 / 
\plot -3.5 -2 -3 -1 /
\plot -2.5 -2 -3 -1 /
\plot -0.5 -1 -1 1 /
\plot -1.5 -1 -1 1 /
\plot 0.5 -1 1 1 /
\plot 1 -1 1 1 /
\plot 1.5 -1 1 1 /

\arrow <0.235cm> [0.2,0.6] from 0 1 to 0.1 1 
\arrow <0.235cm> [0.2,0.6] from  -3 0 to -3 0.1  
\arrow <0.235cm> [0.2,0.6] from -2 0.8 to -2.15 0.815
\arrow <0.235cm> [0.2,0.6] from -2.15 1.2 to -2 1.2 
\arrow <0.235cm> [0.2,0.6] from -3.3 -1.6 to -3.2 -1.4 
\arrow <0.235cm> [0.2,0.6] from -2.7 -1.6 to -2.8 -1.4 
\arrow <0.235cm> [0.2,0.6] from -1.3 -0.2 to -1.25 0   
\arrow <0.235cm> [0.2,0.6] from -0.7 -0.2 to -0.75 0   
\arrow <0.235cm> [0.2,0.6] from 1 -0.235 to 1 0   
\arrow <0.235cm> [0.2,0.6] from 0.7 -0.2 to 0.75 0   
\arrow <0.235cm> [0.2,0.6] from 1.3 -0.2 to 1.25 0   

\setquadratic 
\plot -3 1 -2 0.8 -1 1 /
\plot -3 1 -2 1.2 -1 1 /

\put {$x_1$} [v] at -3 1.2 
\put {$x_2$} [v] at -1 1.2 
\put {$x_3$} [v] at 1 1.2 
\put {$b$} [v] at -3.2 -1 
\put {$\ldots$} at 1.6 -0.1
\endpicture \] 
In $_XE_B$ there are the following extra vertices belonging to $F_E(X,B)$; 
$(f)$ and $(e,f)$ (each emits one edge to $b$), 
$(g)$ and $(e,g)$ (each emits one edge to $x_2$), 
$(h)$, $(d,h)$, $(d,d,h)$,\ldots,$(d,\ldots,d,h)$, \ldots (each 
emits one edge to $x_3$).}
\end{exam} 

We note that, in general, the graph $_XE_B$ has the following structure: 
It contains the restriction $\widetilde{X}=(X,s^{-1}(X),r,s)$ of $E$ 
to $X$. Every vertex in $B\subseteq {_XE_B^0}$ emits infinitely many edges 
into $X$ and does not emit any other edges. Every vertex in $F_E(X,B)$  
emits exactly one edge, which ends in $X\cup B$. Thus every loop in 
the graph $_XE_B$ comes from a loop in $\widetilde{X}$, and 
hence  the loop structure of $_XE_B$ is the same as that of 
$\widetilde{X}$. If $B$ is empty then we write $_XE$ instead of 
$_XE_\emptyset$ and $F_E(X)$ instead of $F_E(X,\emptyset)$. In this 
case we have $F_E(X)=\widetilde{F}_E(X)$. 

As in \cite{bhrs}, for $v\in\Xbad$ we denote by $P_{v,X}$ the projection 
$\displaystyle{\sum_{s(e)=v,r(e)\not\in X}}S_eS_e^*$, and by $J_{X,B}$ 
the ideal of $C^*(E)$ generated by $\{P_w:w\in X\}$ and 
$\{P_v-P_{v,X}:v\in B\}$. According to \cite[Theorem 3.6]{bhrs}, all 
gauge-invariant ideals of $C^*(E)$ are of the form $J_{X,B}$. It has been 
already known that the quotient of $C^*(E)$ by a gauge-invariant ideal is 
itself isomorphic to a graph algebra (cf. 
\cite[Proposition 3.4 and Corollary 3.5]{bhrs}). In the following lemma, 
we show that the same fact holds with respect 
to gauge-invariant ideals. 

\begin{lemm}\label{Jideal} 
Let $E$ be a directed graph. Let $X\neq\emptyset$ be a hereditary and 
saturated subset of $E^0$, and let $B\subseteq\Xbad$. Then there exists 
an isomorphism $$ \phi:C^*(_XE_B)\rightarrow J_{X,B}. $$ 
\end{lemm} 
\begin{proof} 
To avoid confusion, we denote by $T_*$, $Q_*$ the canonical generating family 
of partial isometries and projections for $C^*(_XE_B)$. In order to define 
a $*$-homomorphism $\phi:C^*(_XE_B)\rightarrow J_{X,B}$ 
we first choose target elements for the generating family, as follows.  
For $x\in X$, $b\in B$, $\alpha\in F_E(X,B)$, $e,f\in E^1$ with 
$s(e),r(f)\in X$ and $s(f)\in B$ we set:  
\begin{align*} 
\phi(Q_x) & = P_x \\ 
\phi(Q_b) & = P_b-P_{b,X} \\ 
\phi(Q_{\alpha}) & = \left\{ \hspace{-1.5mm} \begin{array}{ll} 
  S_\alpha S_\alpha^* & \text{ if }\; r(\alpha)\in X \\ S_\alpha(P_{r(\alpha)}-
  P_{r(\alpha),X})S_\alpha^* & \text{ if }\; r(\alpha)\in B \end{array} \right. \\ 
\phi(T_e) & = S_e \\ 
\phi(T_f) & = S_f \\ 
\phi(T_{\overline{\alpha}}) & = \left\{ \hspace{-1.5mm} \begin{array}{ll} 
  S_\alpha & \text{ if }\; r(\alpha)\in X \\ S_\alpha(P_{r(\alpha)}-
  P_{r(\alpha),X}) & \text{ if }\; r(\alpha)\in B \end{array} \right. 
\end{align*} 
It is clear from our definitions that $\{\phi(T_*), \phi(Q_*)\}$ 
forms a Cuntz-Krieger $_XE_B$-family inside $J_{X,B}$. Thus the desired  
$*$-homomorphism exists by the universal property of $C^*(_XE_B)$. 

We show that $\phi$ is injective. To this end we apply the Cuntz-Krieger 
type uniqueness result of \cite[Theorem 1.2]{s2}. Namely, $\phi$ 
is injective if and only if $\phi(Q_v)\neq0$ for each vertex $v\in {_XE_B^0}$, 
and for each vertex-simple loop $\mu$ without exits in $_XE_B$ the spectrum of 
$\phi(T_\mu)$ contains the entire unit circle. (A loop $(e_1,\ldots,e_k)$ 
is vertex-simple if $s(e_i)\neq s(e_j)$ for $i\neq j$.)  
Indeed, the images under $\phi$ 
of all the projections $Q_*$ are non-zero, by construction. 
Furthermore, all loops in $_XE_B$ come from loops in $\widetilde{X}$.  
Thus, if $\mu$ is a vertex-simple loop without exits in $_XE_B$ then 
the spectrum of $\phi(T_\mu)=S_\mu$ contains the entire unit circle 
(cf. \cite[\S2]{kpr}). Consequently $\phi$ is injective. 

We show that $\phi$ is surjective. To this end it sufficies to show that 
every element of the spanning set for $J_{X,B}$, as given in \cite{bhrs}, 
belongs to the range of $\phi$. That is, for $\alpha,\beta$ paths in $E$ 
with $r(\alpha)\in X$ and $r(\beta)\in B$ we must show that $S_\alpha$ 
and $S_\beta(P_{r(\beta)}-P_{r(\beta),X})$ belong to the range of $\phi$. 
Indeed, let $\alpha$ and $\beta$ be such paths in $E$. If $s(\alpha)\in X$ then 
$\alpha$ is a path in $\widetilde{X}$ and $S_\alpha=\phi(T_\alpha)$. If $s(\alpha)
\not\in X$ then there are paths $\alpha_1,\ldots,\alpha_k$ in $E$ such that 
$\alpha=\alpha_1\alpha_2\cdots\alpha_k$, $\alpha_j\in F_E(X,B)$ for 
$j<k-1$, $s(\alpha_k)\in X$ (and possibly $|\alpha_k|=0$), and either 
$\alpha_{k-1}\in F_E(X,B)$ or $\alpha_{k-1}$ is an edge in $E$ from $B$ to 
$X$. In the former case we have $S_\alpha=\phi(T_{\overline{\alpha}_1})\cdots
\phi(T_{\overline{\alpha}_{k-1}})\phi(T_{\alpha_k})$, and in the latter 
$S_\alpha=\phi(T_{\overline{\alpha}_1})\cdots\phi(T_{\overline{\alpha}_{k-2}})
\phi(T_{\alpha_{k-1}})\phi(T_{\alpha_k})$. In all cases $S_\alpha$ belongs 
to the range of $\phi$. Likewise, we have $\beta=\beta_1\cdots\beta_n$ with  
$\beta_j\in F_E(X,B)$, and hence $S_\beta(P_{r(\beta)}-P_{r(\beta),X})=
\phi(T_{\overline{\beta}_1})\cdots\phi(T_{\overline{\beta}_n})$ belongs to 
the range of $\phi$. 
\end{proof} 

The above Lemma \ref{Jideal} nicely complements \cite[Proposition 3.4]{bhrs}. 
Combined, these two results say that the class of graph algebras is closed 
under passing to gauge-invariant ideals and their quotients. It is worth 
noting that according to these two constructions a gauge-invariant ideal $J$ 
of $C^*(E)$ corresponds to a subgraph of $E$ possibly enlarged by some 
extra sources (vertices which receive no edges), while the quotient $C^*(E)/J$ 
corresponds to a subgraph of $E$ possibly enlarged by some extra sinks 
(vertices which emit no edges). 

\section{Type $I$ graph algebras and their limits} 

In this section, our goal is to characterize both graph algebras of type $I$  
and graph algebras which are isomorphic to inductive limits of $C^*$-algebras 
of type $I$. Our characterization of type $I$ graph algebras uses 
the concept of  maximal tail (cf. \cite[Lemma 4.1]{bhrs}). Namely, a 
non-empty subset $M\subseteq E^0$ is a maximal tail if it satisfies 
the following three conditions: 
\begin{description} 
\item[(MT1)] If $v\in E^0$, $w\in M$, and  $v\geq w$, then $v\in M$.  
\item[(MT2)] If $v\in M$ and $0<|s^{-1}(v)|<\infty$, then there 
exists $e\in E^1$ with $s(e)=v$ and $r(e)\in M$. 
\item[(MT3)] For every $v,w\in M$ there exists $y\in M$ such that 
$v\geq y$ and $w\geq y$. 
\end{description}  
(As usual, we write $v\geq w$ if there is a path from $v$ to $w$.) 
The collection of all maximal tails $M$ in $E$ such that each loop in 
$\widetilde{M}=(M,s^{-1}(M),r,s)$ has an exit into $M$ is denoted 
by ${\mathcal M}_\gamma(E)$. That is, if $\alpha=(a_1,\ldots,a_k)$ 
is a loop in $\widetilde{M}$ then there exists an index $j$ and an 
edge $e$ such that $e\neq a_j$, $s(e)=s(a_j)$ and $r(e)\in M$. 
The set of maximal tails $M$ containing a loop without 
exits into $M$ is denoted by ${\mathcal M}_\tau(E)$. Note that a loop in 
$M\in{\mathcal M}_\tau(E)$ might have exits into $E^0\setminus M$. 
Then ${\mathcal M}(E)={\mathcal M}_\gamma(E)\cup{\mathcal M}_\tau(E)$ 
is the set of all maximal tails. As explained in \cite{hs3}, 
the maximal tails in ${\mathcal M}_\gamma(E)$ give rise to 
gauge-invariant primitive ideals of $C^*(E)$, while those in 
${\mathcal M}_\tau(E)$ to non gauge-invariant ones. 

We recall the following notations, introduced in \cite{bhrs}. If $X$ 
is a subset of $E^0$ then $\Omega(X)$ denotes the set of all those 
vertices $w\in E^0\setminus X$ such that there is no path from $w$ to any 
vertex in $X$. A $v\in E^0$ is called a breaking vertex if $|s^{-1}(v)|
=\infty$ and $0<|s^{-1}(v)\setminus r^{-1}(\Omega(v))|<\infty$. 

\begin{theo}\label{typeI} 
Let $E$ be a directed graph. The $C^*$-algebra $C^*(E)$ is of type $I$ 
if and only if for every maximal tail $M\in{\mathcal M}_\gamma(E)$ one 
of the two holds: either (i) $M$ contains a vertex which emits no 
edges into $M$, or (ii) there is an infinite path $(a_1,a_2,\ldots)$ 
such that $s(a_i)\neq s(a_j)\in M$ for $i\neq j$ and each 
vertex $s(a_i)$ emits only one edge into $M$. 
\end{theo} 
\begin{proof} 
($\Rightarrow$) \hspace{1mm} Let $C^*(E)$ be of type $I$, and 
let $M$ be a maximal tail in ${\mathcal M}_\gamma(E)$. Suppose that 
$M$ satisfies neither condition (i) nor (ii) of Theorem \ref{typeI}. We 
will derive a contradiction by constructing an irreducible representation 
of $C^*(E)$ which does not contain any compact operators. To this end, 
we first show that there is an infinite path $\beta=(b_1,b_2,\ldots)$ 
in $\widetilde{M}$ satisfying the following condition. 
\begin{description} 
\item[(\dag)] For each index $k\in\N$ there exists an index $m>k$ and 
an edge $e$ such that $e\neq b_m$, $s(e)=s(b_m)$, and $r(e)\geq s(b_j)$ 
for some $j\in\N$. 
\end{description} 
We construct such a path $\beta$ by the following inductive process. 
Since $M$ fails condition (i), there is an infinite path $\beta_1=
(b_1^1,b_2^1,\ldots)$ in $\widetilde{M}$. We claim that there exists an 
index $i_1$ such that vertex $s(b_{i_1}^1)$ emits at least two 
edges into $M$. Indeed, if $\beta_1$ contains a loop then such 
a vertex exists since every loop in $\widetilde{M}$ has an exit into $M$, 
by assumption. On the other hand, if $\beta_1$ contains no loop then 
$s(b_k^1)\neq s(b_n^1)$ for $k\neq n$. Thus the claim holds in this 
case too, since $M$ fails condition (ii). Set $b_j:=b_j^1$ for 
$1\leq j\leq i_1$. Let $f_1\neq b_{i_1}$ be an edge such that  
$s(f_1)=s(b_{i_1})$ and $r(f_1)\in M$. Since $M$ satisfies (MT3), 
there are two paths $\zeta_1=(b_{i_1},b_{i_1+1},\ldots,b_{m_1})$ and 
$\mu_1=(f_1,\ldots)$ such that $r(\zeta_1)=r(\mu_1)\in M$. Again, since 
$M$ fails (i), there exists an infinite path $\beta_2=(b_{m_1+1}^2,
b_{m_1+2}^2,\ldots)$ in $\widetilde{M}$ which begins at $r(\zeta_1)=r(b_{m_1})$. 
As before, there exists an index $i_2$ such that vertex $s(b_{i_2}^2)$ 
emits at least two edges into $M$. Set $b_j:=b_j^2$ for $m_1+1\leq j
\leq i_2$, and let $f_2\neq b_{i_2}$ be an edge such that $s(f_2)=
s(b_{i_2})$ and $r(f_2)\in M$. As before, there exist paths 
$\zeta_2=(b_{i_2},b_{i_2+1},\ldots,b_{m_2})$ and $\mu_2=(f_2,\ldots)$ 
such that $r(\zeta_2)=r(\mu_2)\in M$. Proceeding inductively in this 
manner, we obtain an infinite path $\beta=(b_1,b_2,\ldots)$ in 
$\widetilde{M}$ which satisfies condition ($\dag$).  

Let $N:=\{v\in M:v\geq s(b_j) \text{ for some } j\}$. It is 
not difficult to verify that $N$ is a maximal tail in $E$. As shown in 
\cite[\S4]{bhrs}, there exists a surjective homomorphism $\pi:C^*(E)
\rightarrow C^*(\widetilde{N})$. Thus, in order to complete the proof, it 
sufficies to show that $C^*(\widetilde{N})$ has an irreducible 
representation $\rho$ disjoint with the compacts. Then $\rho\circ
\pi$ is the required irreducible representation of $C^*(E)$ 
which does not contain any compact operators. The construction of 
such a representation $\rho$ follows. 

Let $A$ be the collection of all infinite paths in $\widetilde{N}$ 
which are shift-tail equivalent to $\beta$. 
(Recall from \cite[\S2]{kprr} that an infinite path $(e_1,e_2,\ldots)$ 
is shift-tail equivalent to $\beta$ if there exist $k,m\in\N$ such that 
$e_{k+i}=b_{m+i}$ for all $i\in\N$.) Then $A$ is an infinite 
set by ($\dag$).  Let $\H$ be a Hilbert space with an 
orthonormal basis $\{\xi_\alpha:\alpha\in A\}$. For a vertex $v$ in 
$\widetilde{N}$ and an edge $e$ in $\widetilde{N}$ we define a 
projection $Q_v$ and a partial isometry $T_e$ in $\B(\H)$, as follows: 
\begin{align}\label{rho}
Q_v(\xi_\alpha) & = 
  \left\{ \hspace{-1.5mm} \begin{array}{ll} \xi_\alpha 
  & \text{if } v=s(\alpha) \\ 0 & \text{otherwise,} 
  \end{array} \right. \\ 
T_e(\xi_\alpha) & = 
  \left\{ \hspace{-1.5mm} \begin{array}{ll} \xi_{(e,c_1,c_2,\ldots)}  
  & \text{if } r(e)=s(\alpha), \text{ where } 
  \alpha=(c_1,c_2,\ldots)\in A \\ 0 & \text{otherwise.} 
  \end{array} \right. \nonumber 
\end{align} 
One can verify that $\{T_e,Q_v\}$ is a Cuntz-Krieger 
$\widetilde{N}$-family. Thus there is a representation 
$\rho:C^*(\widetilde{N})\rightarrow\B(\H)$ such that $\rho(P_v)=Q_v$ 
and $\rho(S_e)=T_e$ for all $v\in\widetilde{N}^0=N$, 
$e\in\widetilde{N}^1=s^{-1}(N)$. 

We show that $\rho$ is irreducible. Indeed, let $U$ be a bounded operator 
on $\H$ which belongs to the commutator of $\rho(C^*(\widetilde{N}))$. 
Let $\alpha=(c_1,c_2,\ldots)$ be an arbitrary element of $A$. We set 
$\alpha[n]:=(c_1,\ldots,c_n)$. Let $\lambda_\mu$, $\mu\in A$ be the complex 
numbers such that 
$$ U\xi_\alpha=\sum_{\mu\in A}\lambda_\mu\xi_\mu. $$ 
Since $T_{\alpha[n]}T_{\alpha[n]}^*\xi_\alpha=\xi_\alpha$ and $U$ 
commutes with $T_{\alpha[n]}$ and $T_{\alpha[n]}^*$, we have 
$$ U\xi_\alpha=UT_{\alpha[n]}T_{\alpha[n]}^*\xi_\alpha=T_{\alpha[n]}
   T_{\alpha[n]}^*U\xi_\alpha=\sum_{\mu\in A}\lambda_\mu T_{\alpha[n]}
   T_{\alpha[n]}^*\xi_\mu. $$ 
Since this identity holds for all $n\in\N$, we must have $U\xi_\alpha=
\lambda_\alpha\xi_\alpha$. We now show that $\lambda_\alpha=\lambda_\beta$. 
Indeed, since $\alpha$ and $\beta$ are shift-tail equivalent there 
exist $k,m\in\N$ such that $c_{k+i}=b_{m+i}$ for all $i\in\N$. Then 
(with $\alpha[n]=(a_1,\ldots,a_n)$ and $\beta[n]=(b_1,\ldots,b_n)$) we have 
$$ \lambda_\beta\xi_\beta=U\xi_\beta=UT_{\beta[m-1]}T_{\alpha[k-1]}^*
   \xi_\alpha=T_{\beta[m-1]}T_{\alpha[k-1]}^*U\xi_\alpha=\lambda_\alpha
   T_{\beta[m-1]}T_{\alpha[k-1]}^*\xi_\alpha=\lambda_\alpha\xi_\beta. $$ 
Thus $U=\lambda_\beta I$ and, consequently, $\rho$ is irreducible. 

We show that $\rho(C^*(\widetilde{N}))\cap\K(\H)=\{0\}$. To this end it 
sufficies to show that the composition $q\circ\rho$ is faithful on 
$C^*(\widetilde{N})$, where $q:\B(\H)\rightarrow\B(\H)/\K(\H)$ is the 
natural quotient map. This follows from the Cuntz-Krieger Uniqueness 
Theorem. Indeed, since every vertex in $N$ connects to $\beta$ and $\beta$ 
satisfies ($\dag$), it follows that every loop in $\widetilde{N}$ 
has an exit. Since for each $v\in N$ the set $\{\alpha\in A:s(\alpha)=v\}$ 
is infinite, the rank of the projection $Q_v$ 
is infinite, and hence $(q\circ\rho)(P_v)=q(Q_v)\neq0$. Therefore 
$q\circ\rho$ is faithful by \cite[Theorem 1.5]{rs}. 

\vspace{2mm} ($\Leftarrow$) \hspace{1mm} 
Suppose that for each maximal tail $M\in{\mathcal M}_\gamma(E)$ either 
(i) or (ii) holds. We must show that if $\varrho:C^*(E)\rightarrow\B(\H_\varrho)$ 
is an irreducible representation then $\K(\H_\varrho)\subseteq\varrho(C^*(E))$. 
Denote $\J:=\ker(\varrho)$. If $\J$ is not gauge-invariant then $\K(\H_\varrho)
\subseteq\varrho(C^*(E))$ by \cite[Theorem 2.10 and Lemma 2.5]{hs3}. Now  
suppose that $\J$ is gauge-invariant. Since $\varrho$ factors through 
a faithful irreducible representation of $C^*(E)/\J$, the inclusion 
$\K(\H_\varrho)\subseteq\varrho(C^*(E))$ will follow if we prove that 
the quotient $C^*(E)/\J$ contains an ideal isomorphic with an algebra 
of compact operators on some Hilbert space. Since $\J$ is gauge-invariant, 
it follows from \cite[Theorem 3.6 and Corollary 3.5]{bhrs} that 
there exists a directed graph $F$ such that $C^*(E)/\J\cong C^*(F)$. 
Thus by virtue of \cite[Proposition 2.1]{kpr} (or, strictly speaking, by 
an obvious extension of that result to the case of arbitrary graphs), 
it sufficies to find a hereditary subset $X\subseteq F^0$ with  
$C^*(\widetilde{X})$ isomorphic to an algebra of compact operators. 

We know from \cite[Theorem 4.7]{bhrs} that either $\J=J_{\Omega(v),
\Ovbad\setminus\{v\}}$ for some breaking vertex $v\in E^0$, or 
$\J=J_{\Omega(M),\OMbad}$ for some maximal tail $M\in{\mathcal M}_\gamma(E)$.  
If $\J=J_{\Omega(v),\Ovbad\setminus\{v\}}$ then the graph $F$ contains 
a sink $\beta(v)$, by \cite[Corollary 3.5]{bhrs}, and we 
may take $X=\{\beta(v)\}$. Then $C^*(\widetilde{X})\cong\C$. 
If $\J=J_{\Omega(M),\OMbad}$ then $F=\widetilde{M}$, 
again by \cite[Corollary 3.5]{bhrs}. By assumption, $M$ satisfies either 
(i) or (ii). If (i) is satisfied then $\widetilde{M}$ contains a 
sink, call it $w$, and $X=\{w\}$ does the job. Finally suppose that 
(ii) holds, and let $\alpha=(a_1,a_2,\ldots)$ be a 
path in $\widetilde{M}$ with the properties described in (ii). 
In this case we may take $X=\{s(a_j):j\in\N\}$, and then 
$C^*(\widetilde{X})$ is the algebra of compact operators on 
a separable Hilbert space. 
\end{proof} 

We now turn to the study of the class of graph algebras which are 
isomorphic to inductive limits of $C^*$-algebras of type $I$. 
It contains all $AF$ graph algebras and in particular all stable $AF$-algebras 
(whose graphs have no loops at all), all stably finite graph algebras 
(whose graphs may have loops but without exits), a variety of graph 
algebras of type $I$, and all the graph algebras considered in \cite{hs1} 
and \cite{hs2} in the context of certain quantum spaces. In the 
following Theorem \ref{limitoftypeI} we give several characterizations 
of those graph algebras. One of them says that these are precisely 
the algebras corresponding to graphs with isolated loops. According to 
our definition in \S1, if $E$ has isolated loops then any two distinct 
vertex-simple loops pass through disjoint sets of vertices. 
(We call two vertex-simple loops {\em distinct} if one is not 
a cyclic permutation of the other.) 

\begin{theo}\label{limitoftypeI} 
Let $E$ be a directed graph. Then the following conditions are 
equivalent. 
\begin{description} 
\item[(i)] $C^*(E)$ is isomorphic to an inductive limit 
of graph $C^*$-algebras of type $I$. 
\item[(ii)] $C^*(E)$ is isomorphic to an inductive limit 
of $C^*$-algebras of type $I$. 
\item[(iii)] $C^*(E)$ does not contain any properly infinite projections. 
\item[(iv)] None of the projections $P_v$, $v\in E^0$, is properly 
infinite in $C^*(E)$. 
\item[(v)] Graph $E$ has isolated loops. 
\end{description} 
\end{theo} 
\begin{proof} 
(i)$\Rightarrow$(ii) \hspace{1mm} 
This is trivial. 

\vspace{2mm} (ii)$\Rightarrow$(iii) \hspace{1mm} 
This follows from \cite[Corollary 2]{c}. 

\vspace{2mm} (iii)$\Rightarrow$(iv) \hspace{1mm} This is trivial. 

\vspace{2mm} (iv)$\Rightarrow$(v) \hspace{1mm} 
Suppose that there are loops $\alpha=(a_1,\ldots,a_m)$ and $\beta=
(b_1,\ldots,b_n)$ in $E$ such that $s(a_1)=s(b_1)$ but $a_1\neq b_1$. 
Set $v=s(a_1)$. Then $S_\alpha$ and $S_\beta$ are partial isometries 
such that $S_\alpha^*S_\alpha=S_\beta^*S_\beta=P_v$ and $S_\alpha S_\alpha^*+
S_\beta S_\beta^*\leq P_v$. Thus $P_v$ is properly infinite. 

\vspace{2mm} (v)$\Rightarrow$(i) \hspace{1mm} 
This we prove in two steps. At first we show that if a graph $E$ is 
finite and satisfies condition (v) then $C^*(E)$ is of type $I$. Then  
we show that if an infinite graph $E$ satisfies condition (v) then 
$C^*(E)$ can be approximated by $C^*$-algebras of finite graphs 
satisfying (v). 

First suppose that $E$ is finite and satisfies (v). We show that $E$ 
fulfills the condition of Theorem \ref{typeI}. So let $M$ be a 
maximal tail in ${\mathcal M}_\gamma(E)$, and suppose for a moment that 
every vertex of $M$ emits at least one edge into $M$. Since $E$ is finite 
this implies that from each vertex of $M$ there is a path to a loop 
in $\widetilde{M}$. In particular, graph $\widetilde{M}$ contains a 
vertex-simple loop $\alpha_1=(a_1,\ldots,a_{k_1})$.  
Since $M$ belongs to ${\mathcal M}_\gamma(E)$, the loop $\alpha_1$ 
has an exit into $M$. Thus there is an edge $e_1\neq a_j$, $j=1,\ldots,k_1$, 
such that $s(e_1)\in\{s(\alpha_j):j=1,\ldots,k_1\}$ and $r(e_1)\in M$. 
Since $E$ satisfies (v) we have $r(e_1)\not\in\{s(\alpha_j):j=1,\ldots,k_1\}$. 
Then there exists a path from $r(e_1)$ to a vertex-simple loop $\alpha_2$ in 
$\widetilde{M}$. Since $E$ satisfies (v) the loops $\alpha_1$ and $\alpha_2$ 
must have disjoint sets of vertices. Again, the loop $\alpha_2$ has an exit 
into $M$ and we can continue this process, constructing an infinite 
sequence of loops in $\widetilde{M}$ with mutually disjoint sets of 
vertices. This contradicts finiteness of $E$. Therefore we must conclude that 
in each maximal tail $M\in{\mathcal M}_\gamma(E)$ there exists a vertex 
which emits no edges into $M$. Hence $E$ satisfies the condition of Theorem 
\ref{typeI} and consequently $C^*(E)$ is of type $I$. 

Now suppose that $E$ is infinite and satisfies (v). We use the method of 
\cite{rs}, described in \S1 of the present paper, to realize $C^*(E)$ as 
an increasing limit $\displaystyle{\lim_{\longrightarrow}}\,C^*(E_{G_n})$ 
of the $C^*$-algebras of finite graphs $E_{G_n}$. Since 
$E$ satisfies (v), so does each of the finite graphs $E_{G_n}$ by Lemma 
\ref{preserve}. Thus each $C^*(E_{G_n})$ is of type $I$ by the preceding 
argument, and consequently $C^*(E)$ is an increasing limit of graph 
$C^*$-algebras of type $I$. 
\end{proof} 

\section{Stable rank of graph algebras} 

The notion of (topological) stable rank for $C^*$-algebras was introduced by 
Rieffel in his seminal paper \cite{r}. For a unital $C^*$-algebra $\mathbf{A}$ 
the stable rank, denoted $\text{sr}(\mathbf{A})$, is the least integer 
$n$ such that the set of $n$-tuples in $\mathbf{A}^n$ which 
generate $\mathbf{A}$ as a left ideal is dense in $\mathbf{A}^n$. 
If such an integer does not exist then $\text{sr}(\mathbf{A})=\infty$. 
If $\mathbf{A}$ does not have identity then its stable rank is 
defined to be that of its minimal unitization (cf. \cite[Definition 1.4]{r}). 
Our calculation of the stable rank of graph algebras relies on 
numerous results of \cite{r}. 

In this section, we completely determine the stable rank of the 
$C^*$-algebra of an arbitrary countable directed graph. 
If $E$ is a directed graph then, by virtue of 
\cite[Proposition 5.5]{rs}, $C^*(E)$ has stable rank one if 
and only if no loop in $E$ has an exit. On the other hand, if $C^*(E)$ 
has a unital, purely infinite and simple quotient, then the stable rank of 
$C^*(E)$ is infinite by \cite[Proposition 6.5 and Theorem 4.3]{r}. Our 
aim is to show that for all other graphs $E$ the stable rank of $C^*(E)$ 
is two. We begin by giving the condition on a graph which guarantees the 
existence of a unital, purely infinite and simple quotient. 

\begin{prop}\label{piquotient} 
Let $E$ be a directed graph. Then there exists an ideal $\J$ of $C^*(E)$ 
with $C^*(E)/\J$ unital, purely infinite and simple if and only if there 
exists a finite $M\in{\mathcal M}_\gamma(E)$ such that $\widetilde{M}$ 
contains at least one loop and does not admit any 
non-trivial hereditary and saturated subsets. 
\end{prop} 
\begin{proof} 
Let $\J$ be an ideal of $C^*(E)$ such that $C^*(E)/\J$ is unital, purely 
infinite and simple. Then $\J$ is primitive and has one of the three forms 
described in \cite[Corollary 2.11]{hs3}. However, it is not possible that 
$\J=J_{\Omega(v),\Ovbad\setminus\{v\}}$ with $v$ a breaking vertex or that 
$\J=R_{N,t}$ with $N\in{\mathcal M}_\tau(E)$ and $t\in\T$. Indeed, The 
quotient $C^*(E)/R_{N,t}$ contains an ideal isomorphic with the compacts by 
\cite[Lemma 2.5]{hs3}. Likewise, the graph corresponding to $C^*(E)/
J_{\Omega(v),\Ovbad\setminus\{v\}}$ contains a sink by \cite[Corollary 3.5]{bhrs} 
and hence this quotient also contains an ideal isomorphic with the compacts. 
Thus we must have  $\J=J_{\Omega(M),\OMbad}$ for some $M\in{\mathcal M}_\gamma(E)$. 
In this case $C^*(E)/\J$ is isomorphic to $C^*(\widetilde{M})$. So we need 
the necessary and sufficient conditions on $\widetilde{M}$ with 
$M\in{\mathcal M}_\gamma(E)$ which guarantee that $C^*(\widetilde{M})$ 
is unital, purely infinite and simple. 

Clearly, $C^*(\widetilde{M})$ is 
unital if and only if $M$ is finite. Since $M\in{\mathcal M}_\gamma(E)$, 
every loop in $\widetilde{M}$ has an exit. Thus $C^*(\widetilde{M})$ is 
simple if and only if $\widetilde{M}$ does not admit any non-trivial 
hereditary and saturated sets, by \cite[Theorem 12]{s1}. Furthermore, 
$C^*(\widetilde{M})$ is infinite if and only if $\widetilde{M}$ contains 
at least one loop  with an exit. If $C^*(\widetilde{M})$ is simple 
this implies that $C^*(\widetilde{M})$ is purely infinite 
(cf. \cite[Corollary 3.11]{kpr} and \cite[Theorem 18]{s1}). 
\end{proof} 

Since the stable rank of any stable $C^*$-algebra is either one or two 
by \cite[Theorem 6.4]{r}, the following Lemma \ref{stableideal} will be 
a key ingredient in our determination of stable rank of graph algebras. 
The lemma says that in absence of unital, purely infinite and simple quotients 
one can find a stable ideal $\J$ of $C^*(E)$ such that $C^*(E)/\J$ does not 
contain any properly infinite projections.  

\begin{lemm}\label{stableideal} 
Let $E$ be a directed graph. If $C^*(E)$ does not have any unital, purely 
infinite and simple quotients, then there exists a stable, gauge-invariant 
ideal $\J$ of $C^*(E)$ such that the quotient $C^*(E)/\J$ is isomorphic 
to the algebra of a graph with isolated loops. 
\end{lemm} 
\begin{proof} 
Let $E$ be a directed graph such that $C^*(E)$ does not have any unital, purely 
infinite and simple quotients. We denote by $X_0$ the collection of all those 
vertices $v\in E^0$ for which there exist two edges $e\neq f\in E^1$ such 
that $s(e)=s(f)=v$, $r(e)\geq v$ and $r(f)\geq v$. We let $X$ be the smallest 
hereditary and saturated subset of $E^0$ containing $X_0$. We show that the 
ideal $\J=I_X$ generated by the projections $P_v$, $v\in X$, has the desired 
properties. Clearly, $\J$ is invariant under the gauge action. By 
\cite[Proposition 3.4]{bhrs} the quotient $C^(E)/\J$ is isomorphic to $C^*(E/X)$. 
The graph $E/X$ is the restriction of $E$ to $E^0\setminus X$, with a 
possible addition of some extra sinks (cf. \cite{bhrs}). It follows that 
$E/X$ has isolated loops. Indeed, all loops of $E/X$ pass through  
vertices in $E^0\setminus X$. Thus, if they were not isolated, the 
intersection of $E^0\setminus X$ and $X_0$ would not be empty, a contradiction. 

It remains to show that the ideal $\J$ is stable. 
By Lemma \ref{Jideal}, $\J$ is isomorphic to the graph algebra 
$C^*(_XE)$. Hence, by virtue of \cite[Theorem 3.2]{t}, in order 
to show that $\J$ is stable it sufficies to prove that 
for every loop $\alpha$ in $_XE$ there are infinitely many vertices 
$w\in{_XE}$ connecting to $\alpha$, and $_XE$ has no non-zero bounded 
graph traces. We prove both these properties by contradiction. 

Suppose that there is a loop $\alpha=(a_1,\ldots,a_k)$ in $_XE$ such that 
there are only finitely many vertices in $_XE$ connecting to 
$\{r(a_j):j=1,\ldots,k\}$. Let $\Omega(\alpha)$ be the set of all 
those vertices in $_XE$ from which there is no path to $\{r(a_j):j=1,\ldots,k\}$. 
Then $_XE^0\setminus\Omega(\alpha)$ is a finite maximal tail in $_XE$. 
Let $M$ be a maximal tail of smallest cardinality contained in $_XE^0
\setminus\Omega(\alpha)$. Since $M$ is finite, $C^*(\widetilde{M})$ 
is unital. We have $M\cap X_0\neq\emptyset$, for otherwise 
$X\setminus M$ would be a hereditary and saturated proper subset of $X$ 
containing $X_0$, a contradiction. So let $v\in M\cap X_0$. If $w\in M$ 
lies on a loop in $\widetilde{M}$ then $w\geq v$. For otherwise 
$M\setminus\Omega(v)$ would be a maximal tail in $_XE^0\setminus\Omega(\alpha)$ 
with fewer elements than $M$. Hence all loops in $\widetilde{M}$ 
have exits. Since $M$ does not contain any smaller maximal tails, 
it follows from \cite[Corollary 2.11]{hs3} 
that $C^*(\widetilde{M})$ has no non-zero primitive ideals. Therefore 
$C^*(\widetilde{M})$ is simple. Since $M$ contains an element $v$ of $X_0$, 
the graph $\widetilde{M}$ contains the loops passing through $v$. 
This implies that $C^*(\widetilde{M})$ is purely infinite 
(cf. \cite[Corollary 3.11]{kpr} and \cite[Theorem 18]{s1}). 
Hence $\J\cong C^*(_XE)$ has a unital, purely 
infinite and simple quotient $C^*(_XE)/J_{\Omega(M),\OMbad}\cong 
C^*(\widetilde{M})$. Consequently, $C^*(E)$ has such a quotient as well, 
a contradiction. 

Suppose that $\psi$ is a bounded graph trace on $_XE$. According to 
\cite[Definition 2.2]{t}, this means that $\psi: {_XE^0}\rightarrow\R^+$ is 
a function such that  
\begin{description} 
\item[(GT1)] $\psi(v)=\displaystyle{\sum_{s(e)=v}}\psi(r(e))$ for 
all $v\in {_XE^0}$ such that $0<|s^{-1}(v)|<\infty$,  
\item[(GT2)] $\psi(v)\geq\displaystyle{\sum_{s(e)=v}}\psi(r(e))$ for 
all $v\in {_XE^0}$ such that $|s^{-1}(v)|=\infty$,  
\end{description} 
and $\displaystyle{\sum_{v\in {_XE^0}}}\psi(v)<\infty$. 
It follows from (GT1) and (GT2) that $\psi(v)=0$ for every $v\in X_0$. 
Indeed, let $\alpha=(a_1,\ldots,a_m)$ and $\alpha=(b_1,\ldots,b_n)$ be two 
loops in $_XE$ such that $s(\alpha)=s(\beta)=v$ but $a_1\neq b_1$. Then $\psi(v)
\geq\psi(r(a_1))\geq\ldots\geq\psi(r(a_m))=\psi(v)$, and hence $\psi(v)=
\psi(r(a_1))$. Likewise $\psi(v)=\psi(r(b_1))$. But then $\psi(v)\geq
\psi(r(a_1))+\psi(r(b_1))=2\psi(v)$ and hence $\psi(v)=0$, as claimed. 
Consequently, $\psi(w)=0$ for all $w\in {_XE^0}$ by \cite[Lemma 3.7]{t}, 
since the smallest hereditary and saturated subset of $_XE^0$ containing 
$X_0$ is the entire $_XE^0$. 
\end{proof} 

The following example illustrates the construction of the ideal $\J$ 
in the preceding Lemma \ref{stableideal}, and gives the graphs for 
$\J$ and the corresponding quotient. 

\begin{exam}\label{example2} 
{\rm Let $E$ be the following graph: 

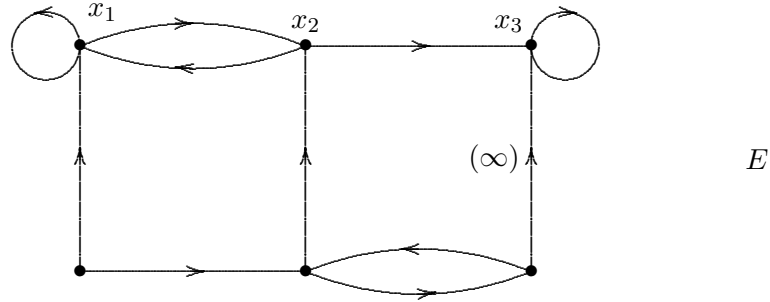
\[ \beginpicture
\setcoordinatesystem units <1.5cm,1.5cm>
\setplotarea x from -4 to 1, y from -1.6 to 1 
\put {$E$} at 3 0

\put {$\bullet$} at -3 -1
\put {$\bullet$} at -3 1
\put {$\bullet$} at -1 -1
\put {$\bullet$} at -1 1
\put {$\bullet$} at 1 -1
\put {$\bullet$} at 1 1

\setlinear  
\plot -3 -1 -3 1 /
\plot -1 -1 -1 1 /
\plot -1 1 1 1 / 
\plot 1 -1 1 1 /
\plot -3 -1 -1 -1 / 

\arrow <0.235cm> [0.2,0.6] from 0 1 to 0.1 1 
\arrow <0.235cm> [0.2,0.6] from 1 0 to 1 0.1 
\arrow <0.235cm> [0.2,0.6] from  -3 0 to -3 0.1  
\arrow <0.235cm> [0.2,0.6] from -1 0 to -1 0.1 
\arrow <0.235cm> [0.2,0.6] from -2 -1 to -1.9 -1 

\setquadratic 
\plot -3 1  -2 0.8  -1 1 /
\plot -3 1  -2 1.2  -1 1 / 
\plot -1 -1  0 -1.2  1 -1 / 
\plot -1 -1  0 -0.8  1 -1 / 

\arrow <0.235cm> [0.2,0.6] from -2 0.8 to -2.15 0.815
\arrow <0.235cm> [0.2,0.6] from -2.15 1.2 to -2 1.2 

\arrow <0.235cm> [0.2,0.6] from 0 -1.2 to 0.15 -1.185 
\arrow <0.235cm> [0.2,0.6] from 0 -0.8 to -0.15 -0.8  

\circulararc 360 degrees from 1 1 center at 1.3 1  
\arrow <0.236cm> [0.2,0.6] from 1.3 1.3 to 1.4 1.295 

\circulararc 360 degrees from -3 1 center at -3.3 1  
\arrow <0.236cm> [0.2,0.6] from -3.3 1.3 to -3.4 1.295

\put {$x_1$} [v] at -2.8 1.3 
\put {$x_2$} [v] at -1 1.2 
\put {$x_3$} [v] at 0.8 1.2 
\put {$(\infty)$} [r] at 0.9 0 

\endpicture \] 
One can check with help of Proposition \ref{piquotient} that 
$C^*(E)$does not have any unital, purely infinite and simple quotients. 
We have $X_0=\{x_1\}$ and $X=\{x_1, x_2, x_3\}$. The ideal $\J=I_X$ of 
$C^*(E)$ is stable and, by Lemma \ref{Jideal}, isomorphic to the $C^*$-algebra of 
the following graph $_XE$ (in which both $x_2$ and $x_3$ receive infinitely 
many edges, each beginning at a different vertex): 

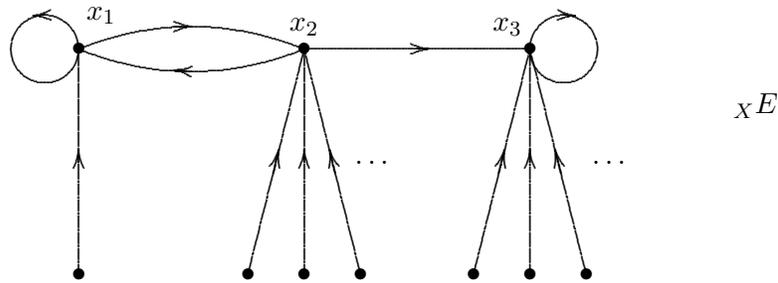
\[ \beginpicture
\setcoordinatesystem units <1.5cm,1.5cm>
\setplotarea x from -4 to 1, y from -1.2 to 1.5 

\put {$_XE$} at 3 0.5 

\put {$x_1$} [v] at -2.8 1.3 
\put {$x_2$} [v] at -1 1.2 
\put {$x_3$} [v] at 0.8 1.2 

\put {$\bullet$} at -3 -1
\put {$\bullet$} at -3 1
\put {$\bullet$} at -1 -1
\put {$\bullet$} at -1.5 -1 
\put {$\bullet$} at -0.5 -1
\put {$\bullet$} at -1 1
\put {$\bullet$} at 1 -1
\put {$\bullet$} at 1.5 -1
\put {$\bullet$} at 0.5 -1 
\put {$\bullet$} at 1 1 

\put {$\ldots$} at 1.7 0
\put{$\ldots$} at -0.4 0 

\setlinear  
\plot -3 -1 -3 1 /
\plot -1 -1 -1 1 /
\plot -1.5 -1 -1 1 /
\plot -0.5 -1 -1 1 /
\plot -1 1 1 1 / 
\plot 1 -1 1 1 /
\plot 1.5 -1 1 1 /
\plot 0.5 -1 1 1 / 

\arrow <0.235cm> [0.2,0.6] from  -3 0 to -3 0.1  
\arrow <0.235cm> [0.2,0.6] from 0 1 to 0.1 1 
\arrow <0.235cm> [0.2,0.6] from -1 0 to -1 0.1 
\arrow <0.235cm> [0.2,0.6] from 1 0 to 1 0.1 
\arrow <0.235cm> [0.2,0.6] from -1.275 -0.1 to -1.225 0.1 
\arrow <0.235cm> [0.2,0.6] from 0.725 -0.1 to 0.775 0.1 
\arrow <0.235cm> [0.2,0.6] from -0.725 -0.1 to -0.775 0.1 
\arrow <0.235cm> [0.2,0.6] from 1.275 -0.1 to 1.225 0.1 

\setquadratic 
\plot -3 1 -2 0.8 -1 1 /
\plot -3 1 -2 1.2 -1 1 /

\arrow <0.235cm> [0.2,0.6] from -2 0.8 to -2.15 0.815
\arrow <0.235cm> [0.2,0.6] from -2.15 1.2 to -2 1.2  

\circulararc 360 degrees from 1 1 center at 1.3 1  
\arrow <0.236cm> [0.2,0.6] from 1.3 1.3 to 1.4 1.295  

\circulararc 360 degrees from -3 1 center at -3.3 1  
\arrow <0.236cm> [0.2,0.6] from -3.3 1.3 to -3.4 1.295 

\endpicture \] 
The quotient $C^*(E)/\J$ is isomorphic to the $C^*$-algebra 
of the following graph $E/X$: 
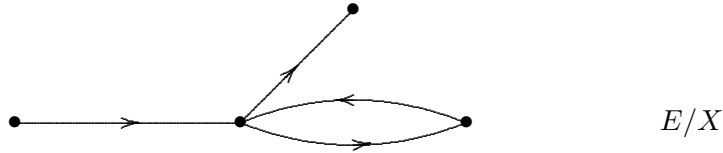
\[ \beginpicture
\setcoordinatesystem units <1.5cm,1.5cm>
\setplotarea x from -4 to 1, y from -0.6 to 1 

\put {$E/X$} at 3 0

\put {$\bullet$} at -3 0
\put {$\bullet$} at -1 0
\put {$\bullet$} at 1 0 
\put {$\bullet$} at 0 1 

\setlinear  
\plot -3 0  -1 0  0 1 /

\arrow <0.235cm> [0.2,0.6] from -2 0 to -1.9 0 
\arrow <0.235cm> [0.2,0.6] from -0.6 0.4 to -0.5 0.5   

\arrow <0.235cm> [0.2,0.6] from 0 -0.2 to 0.15 -0.185 
\arrow <0.235cm> [0.2,0.6] from 0 0.2 to -0.15 0.2  

\setquadratic 
\plot -1 0  0 -0.2  1 0 / 
\plot -1 0  0 0.2  1 0 / 

\endpicture \] 
}
\end{exam} 

Now we are ready to prove our main result. 

\begin{theo}\label{main} 
Let $E$ be a directed graph. Then 
$$ \text{sr}(C^*(E))=\left\{ \hspace{-1.5mm} \begin{array}{ll} 
  1 & \text{ if }\; \text{ no loop in $E$ has an exit}, \\ 
  \infty & \text{ if }\; \text{ $C^*(E)$ has a unital, purely infinite 
  and simple quotient}, \\ 
  2 & \text{ otherwise.} \end{array} \right. $$ 
\end{theo} 
\begin{proof} 
For a directed graph $E$, the stable rank of $C^*(E)$ is one if and only if 
no loop in $E$ has an exit by \cite[Proposition 5.5]{rs}. On the other hand, 
if $C^*(E)$ has a unital, purely infinite and simple quotient then 
its stable rank is infinite by \cite[Theorem 4.3 and Proposition 6.5]{r}. 
Therefore, it remains to prove that $\text{sr}(C^*(E))\leq2$ if $C^*(E)$ 
does not admit any unital, purely infinite and simple quotients. 
So assume this is the case, and let $\J$ be a stable, gauge-invariant 
ideal of $C^*(E)$ such that the quotient $C^*(E)/\J$ is isomorphic to 
the $C^*$-algebra of a graph $D$ with isolated loops. Such an ideal $\J$ 
exists by Lemma \ref{stableideal}. We have $\text{sr}(\J)\leq2$ 
by \cite[Theorem 6.4]{r}. 

At first we consider the case when $D$ has only finitely many distinct vertex-simple 
loops, and proceed by induction on the number $L$ of such loops. Namely, we 
prove that if $\mathbf{A}$ is a $C^*$-algebra which admits an exact sequence 
\begin{equation}\label{jadextension} 
0 \longrightarrow \J \longrightarrow \mathbf{A} 
\stackrel{\pi}{\longrightarrow} C^*(D) \longrightarrow 0, 
\end{equation}  
where $\text{sr}(\J)\leq2$ and $D$ is a directed graph with finitely many 
isolated loops, then $\text{sr}(\mathbf{A})\leq2$. If $L=0$ then 
$C^*(D)$ is $AF$ (cf. \cite[Theorem 2.4]{kpr} and \cite[\S5.4]{rs}) and has 
stable rank $1$ (cf. \cite[Proposition 3.5]{r}). Thus in the exact sequence 
(\ref{jadextension}) we have $\text{sr}(\J)\leq2$ and $\text{sr}(C^*(D))=1$, 
and consequently $\text{sr}(\mathbf{A})\leq\max(\text{sr}(\J),\text{sr}
(C^*(D))+1)=2$ by \cite[Corollary 4.12]{r}. 

For the inductive step, we consider the binary relation $\geq$ defined on the set 
of distinct vertex-simple loops of $D$ by $\alpha\geq\beta$ if there is a 
path from $\alpha$ to $\beta$. This relation is a partial order since the 
loops in $D$ are isolated. Since $L$ is finite there exists a loop $\alpha$ 
which is maximal with respect to $\geq$. Consider the ideal $J_{\Omega(\alpha),
\Oabad}$ of $C^*(D)$. By \cite[Corollary 3.5]{bhrs}, the quotient 
$C^*(D)/J_{\Omega(\alpha),\Oabad}$ is isomorphic to $C^*(D\setminus\Omega
(\alpha))$, where $D\setminus\Omega(\alpha)$ denotes the restriction 
of graph $D$ to $D^0\setminus\Omega(\alpha)$. By Lemma \ref{Jideal}, 
we have $J_{\Omega(\alpha),\Oabad}\cong C^*( {_{\Omega(\alpha)}D_{\Oabad}})$. 
Thus (\ref{jadextension}) implies existence of the exact sequences  
\begin{equation}\label{ext1} 
0 \longrightarrow \pi^{-1}(J_{\Omega(\alpha),\Oabad}) \longrightarrow 
\mathbf{A} \longrightarrow C^*(D\setminus\Omega(\alpha)) \longrightarrow 0  
\end{equation} 
and 
\begin{equation}\label{ext2} 
0 \longrightarrow \J \longrightarrow \pi^{-1}(J_{\Omega(\alpha),\Oabad}) 
\longrightarrow C^*( {_{\Omega(\alpha)}D_{\Oabad}}) \longrightarrow 0. 
\end{equation} 
By our construction, the graph $_{\Omega(\alpha)}D_{\Oabad}$ has one 
distinct vertex-simple loop less than $D$ (cf. Definition \ref{idealgraph}). 
Hence $\text{sr}(\pi^{-1}(J_{\Omega(\alpha),\Oabad}))\leq2$ by the 
inductive hypothesis. Also, by our choice of $\alpha$, the graph 
$D\setminus\Omega(\alpha)$ has a unique (up to a cyclic permutation) 
vertex-simple loop, namely $\alpha$, and this loop has no exits in 
$D\setminus\Omega(\alpha)$. Thus $\text{sr}(C^*(D\setminus\Omega(\alpha)))
=1$ by \cite[Proposition 5.5]{rs}. Therefore the exact sequence 
(\ref{ext1}) implies that $\text{sr}(\mathbf{A})\leq2$, and the inductive 
step follows. 

Now we consider the general case, with possibly infinitely many distinct 
vertex-simple loops in the graph $D$. We have the exact sequence 
$$ 0 \longrightarrow \J \longrightarrow C^*(E)  
   \stackrel{\pi}{\longrightarrow} C^*(D) \longrightarrow 0. $$  
If $G$ is a finite subset of $D^0\cup D^1$ then the graph $D_G$ has isolated 
loops by Lemma \ref{preserve}. Since $D_G$ is finite, so is the number 
of distinct vertex-simple loops in it. By Lemma \ref{approximation}, 
there is an exact sequence 
$$ 0 \longrightarrow \J \longrightarrow \pi^{-1}(C^*(D_G))  
   \longrightarrow C^*(D_G) \longrightarrow 0. $$ 
Thus $\text{sr}(\pi^{-1}(C^*(D_G)))\leq2$ by the preceding argument. As shown 
in \S1 (cf. Lemma \ref{approximation} and the comments following it), 
$C^*(D)$ is isomorphic to the increasing limit of algebras of the form 
$C^*(D_G)$. Thus $C^*(E)$ is the increasing limit of the algebras 
$\pi^{-1}(C^*(D_G))$. Consequently, we have $\text{sr}(C^*(E))
\leq\lim\inf(\text{sr}(\pi^{-1}(C^*(D_G))))=2$ by \cite[Theorem 5.1]{r}. 
Finally, if $E$ has a loop with an exit then $C^*(E)$ is infinite, 
and hence $\text{sr}(C^*(E))\neq1$. Thus in this case $\text{sr}(C^*(E))=2$. 
If there is no loop in $E$ with an exit then $\text{sr}(C^*(E))=1$ by 
\cite[Proposition 5.5]{rs}. 
\end{proof}

\end{document}